\documentclass[a4paper,11pt]{article}
\usepackage{latexsym}
\usepackage{amssymb}
\usepackage{amsfonts}
\usepackage{amsmath}
\usepackage{indentfirst}
\usepackage{graphicx}
\usepackage{epsfig}
\usepackage{color}
\usepackage{fullpage}
\usepackage{tikz,ifthen}

\newtheorem{prop}{Proposition}[section]
\newtheorem{teor}{Theorem}[section]

\newtheorem{cor}{Corollary}[section]

\newcommand{\cvd}{\hfill $\blacksquare$\bigskip}
\date{}

\author{Luca Ferrari\thanks{Dipartimento di Matematica e Informatica ``U. Dini", viale Morgagni 65, University
of Firenze, Firenze, Italy, {\tt luca.ferrari@unifi.it}. Member of the INdAM Research group GNCS;
partially supported by the 2019 INdAM-GNCS project ``Studio di propriet\'a combinatoriche di linguaggi formali ispirate dalla biologia e da
strutture bidimensionali'' and by a grant of the ``Fondazione della Cassa di Risparmio di Firenze'' for the project
``Rilevamento di pattern: applicazioni a memorizzazione basata sul DNA, evoluzione del genoma, scelta sociale''.}}

\title{Enhancing the connections between patterns in permutations and forbidden configurations in restricted elections}

\begin{document}

\maketitle

\begin{abstract}
We investigate the connections between patterns in permutations and forbidden configurations in restricted elections,
first discovered by Lackner and Lackner, in order to enhance the approach initiated by the two mentioned authors.
More specifically, our achievements are essentially two.
First, we define a new type of domain restriction, called \emph{enriched group-separable}.
Enriched group-separable elections are a subset of group-separable elections,
which describe a special, still natural, situation that can arise in the context of group-separability.
The exact enumeration of group-separable elections has been very recently determined by Karpov.
Here we give a recursive characterization for enriched group-separable elections,
from which we are able to find a recurrence relation and a closed formula expressing their number.
Our second achievement is a generalization of a result of Lackner and Lackner,
concerning the connection between permutation patterns and forbidden configurations with 3 voters.
Our result relates forbidden configurations with the strong order on pairs of permutations,
a notion which is still largely undeveloped,
and suggests a potential approach for the determination of upper bounds for restricted elections
whose forbidden configurations contains at least one configuration with 3 voters.
\end{abstract}

\section{Introduction}

An \emph{election} is a pair $(C,P)$ where $C$ is a finite set and $P$ is a tuple of total orders on $C$.
The set $C$ is the set of \emph{candidates} and the tuple $P$ is the tuple of \emph{preferences}.
If we are not interested in the specific order in which the preferences are listed, we speak of the \emph{multiset} of the preferences,
where the term ``multiset" refers to the fact the each preference can appear more than once.
When $|C|=m$ and $|P|=n$, we say that $(C,P)$ is a \emph{$(m,n)$-election}.
In this case, the set of candidates is usually the set $[m]=\{ 1,2,\ldots ,m\}$ and the tuple of preferences is denoted $P=(V_1 ,V_2 ,\ldots ,V_n )$.
The set of \emph{electors, or voters} is the set $[n]=\{ 1,2,\ldots ,n\}$, i.e. the set of indices of the preferences.
Each preference $V_i$ represents the ranking of the candidates proposed by voter $i$, and will be sometimes identified with the voter itself.
When voter $i$ prefers candidate $x$ to candidate $y$, we write $x>_{V_i}y$.
We explicitly notice that, throughout the whole paper, rankings are assumed to be total orders, and in particular they do not contain ties.

\bigskip

In social choice theory it is rather customary to apply some \emph{domain restriction} to elections.
This is mainly motivated by the attempt to find reasonably interesting elections that escape certain classical paradoxes
(such as Arrow's Impossibility Theorem or Condorcet paradox).
Well known and studied domain restrictions are the \emph{single-peaked restriction} \cite{Bl} and the \emph{group-separable restriction} \cite{I}.

It turns out that several domain restrictions can be described in terms of \emph{forbidden configurations}.
Consider a $(m,n)$-election $(C,P)$ and a $(h,l)$-election $(S,T)$, with $T=(T_1 ,T_2 ,\ldots ,T_l )$.
We say that $(C,P)$ \emph{contains the configuration} $(S,T)$ when there exist injective functions $f:[l]\rightarrow [n]$ and $g:S\rightarrow C$ such that
for every $x,y \in S$ and for every $i\leq l$, if $x>_{T_i}y$ then $g(x)>_{V_{f(i)}}g(y)$.
Otherwise, we say that $(C,P)$ \emph{avoids} $(S,T)$.

In recent years there have been quite a few publications dealing with domain restrictions and forbidden configurations.
Both single-peaked and group-separable domains can be characterized in terms of forbidden configurations \cite{BH}.
In \cite{BCW} a similar characterization has been provided for single-crossing domains,
whereas one-dimensional Euclidean domains are studied in \cite{CG,CPW}.
Also, some investigation concerning the enumeration of elections with given numbers of voters and candidates have been pursued.
For instance, the number of certain single-peaked and single-crossing election has been determined in \cite{CF}.
Very recently, a closed form for the number of group-separable elections has been obtained in \cite{K},
by using combinatorial tools like parenthesizations and lattice paths (in particular Schr\"oder paths).

In Section \ref{enriched} we will explicitly state the forbidden configuration characterization of group-separable elections found in \cite{BH}.
More recently,
Lackner and Lackner \cite{LL} have revealed the intimate connection between
forbidden configurations in elections and \emph{pattern avoidance} in permutations
and have exploited it to find interesting combinatorial and probabilistic results on the above mentioned domain restrictions.
In the next section we provide a quick overview of those achievements of \cite{LL} that are relevant for our paper.
Below we give the main definitions concerning patterns in permutations.
Let $S_n$ be the symmetric group over a set of cardinality $n$, consisting of all permutations of length $n$.
Given two permutations $\sigma \in S_k$ and $\tau =\tau_1 \tau_2 \cdots \tau_n \in S_n$, with $k\leq n$,
we say that $\sigma$ is a \emph{pattern} of $\tau$ when
there exist $1\leq i_1 <i_2 <\cdots <i_k \leq n$ such that
$\tau_{i_1}\tau_{i_2}\cdots \tau_{i_k}$ (as a permutation) is isomorphic to $\sigma$
(which means that $\tau_{i_1}, \tau_{i_2}, \ldots , \tau_{i_k}$ are in the same relative order
as the elements of $\sigma$).
When $\sigma$ is not a pattern of $\tau$, we say that $\tau$ \emph{avoids} $\sigma$.
Thus, for instance, 312 is a pattern of $\tau =526143$ (as shown, for example, by the three elements 5,1 and 3);
however, $\tau$ avoids 123.

\bigskip

The present paper builds on \cite{LL}
and illustrates a couple of results that strengthen the validity of the approach initiated in the mentioned paper.

Our first achievement is the definition of a new kind of domain restriction, called \emph{enriched group-separable},
which adds more constraints to the group-separable domain.
As such, an enriched group-separable election is even less likely to appear than a group-separable one,
yet we claim that it has a nontrivial meaning from the point of view of social choice.
Moreover, its links with pattern avoidance in permutations allows us to find nice exact enumerative results
that are much harder to achieve for general group-separable elections.

The second result of our paper is a slight generalization of a result of \cite{LL} concerning $(m,3)$-elections
to the case in which all the preferences are different.
It turns out that, in such a situation, configuration containment can be equivalently expressed in terms of
\emph{strong containment of pairs of permutations}.
This notion for permutations is not new, but it has been little studied, and very few results are known
(probably because of its intrinsic difficulty).
We hope that the link with configuration containment we have found can stimulate further work on this difficult partial order on pairs of permutations.

\section{Preliminaries}

The goal of this section is to briefly describe the connections between domain restrictions and permutation patterns that have been investigated in \cite{LL}.
All the stated results come from such a paper, so we avoid to explicitly recall it for each of them.
In this context, the main result is the following, which shows that the desired link involves elections with 3 voters.
From now on, we use the notation $\tau \leq \pi$ to mean that $\tau$ is a pattern of $\pi$.

\begin{prop}\label{LandL}
Let $\tau =\tau_1 \tau_2 \cdots \tau_h \in S_h$ and $\pi =\pi_1 \pi_2 \cdots \pi_m \in S_m$, with $h\leq m$.
Suppose that $(C,P)$ is a $(m,3)$-election, with $C=[m]=\{ 1,2,\ldots m\}$ and $P=(V_1 ,V_2 ,V_3 )$,
such that $V_1 :12\cdots m$, $V_2 :12\cdots m$ and $V_3 =\pi_1 \pi_2 \cdots \pi_m$
(here, as usual in this context, total orders are represented by listing elements in decreasing order).
Moreover, let $(S,T)$ be a $(h,3)$-election, with $S=[h]$ and $T=(T_1 ,T_2 ,T_3 )$,
such that $T_1 :12\cdots h$, $T_2 :12\cdots h$ and $T_3 =\tau_1 \tau_2 \cdots \tau_h$.
Then $(C,P)$ contains the configuration $(S,T)$ if and only if $\tau \leq \pi$.
\end{prop}

As a consequence, the following proposition finds an analogous characterization for elections with 2 voters.

\begin{prop}
Let $(C,P)$ and $(S,T)$ be elections as in the previous proposition, but with 2 voters.
More specifically, in the notations above just remove $V_2$ and $T_2$ and rename $V_3$ and $T_3$ as $V_2$ and $T_2$, respectively.
Then $(C,P)$ contains the configuration $(S,T)$ if and only if $\tau \leq \pi$ or $\tau^{-1}\leq \pi$,
where $\tau^{-1}$ denotes the group-theoretic inverse of the permutation $\tau$.
\end{prop}

The importance of the last proposition lies basically in the fact that it allows to count restricted elections, at least in the case of 2 voters.
Just a piece of notation: for a given set of permutations $\Pi$,
the set of all permutations of length $n$ avoiding each permutation of $\Pi$ is denoted $S_n (\Pi )$.

\begin{cor}\label{inverse}
Let $(S,T)$ be a $(h,2)$-election, with $T=(T_1 ,T_2 )$.
Let $C=[m]$ and $V_1$ be a total order on $C$.
Then the number of total orders $V_2$ on $C$ such that the election $(C,\{ V_1 ,V_2 \} )$ avoids $(S,T)$ is equal to $|S_n (\tau ,\tau^{-1})|$,
where $\tau$ is the permutation which maps $T_1$ into $T_2$.
\end{cor}

For more general elections (i.e., with more than 2 voters),
the above results can be used to find nontrivial upper bounds for the number of restricted elections.

\begin{teor}
Denote with $a(n,m,\Gamma )$ the number of all $(m,n)$-elections avoiding the set of configurations $\Gamma$.
If $\Gamma$ contains a $(k,2)$-configuration, for $k\geq 2$, then, for all $n,m \in \mathbf{N}$,
$$
a(n,m,\Gamma )\leq m!\cdot c_k ^{(n-1)m},
$$
where $c_k$ only depends on $k$.
\end{teor}

The above theorem allows Lackner and Lackner to find interesting upper bounds, for instance, for the number of single-peaked and group-separable elections,
which we do not report here for the sake of brevity.

\bigskip

We close this section with two remarks.

It is implicit in Corollary \ref{inverse} that, given a set of permutations that is closed under inversion,
we can translate these permutations into configurations (consisting of two conditions).
If we take ``interesting" sets of permutations (closed under inversion), do these sets correspond to reasonable domain restrictions?
Namely, what kind of restrictions does such a set of permutations impose on a set of total orders that avoid the corresponding configurations?
The content of the next section proposes an instance of this general problem that we believe to be quite interesting.

It is often the case that configurations can be conveniently represented using posets.
More precisely, a set of configuration can often be interpreted as the set of all pairs of linear extensions of suitable posets.
Rather than giving a more formal explanation of this statement,
we refer to Section \ref{enriched} for an explicit description of this kind of enriched group-separable elections.

\section{Enriched group-separable elections}\label{enriched}

Recall that a \emph{group-separable election} is an election $(C,P)$ such that, for any subset of candidates $C'\subseteq C$,
there exists a partition of $C'$ into two blocks $A$ and $B$ such that, for each elector $i$,
either $i$ ranks all candidates of $A$ above all candidates of $B$ or vice versa.

There is a characterization of group-separable elections in terms of forbidden configurations,
which is due to Ballester and Haeringer \cite{BH}.

\begin{prop}
An election $(C,P)$ is group-separable if and only if the following two conditions hold:
\begin{itemize}
\item for every $x_1 ,x_2 ,x_3 \in C$, there cannot be three distinct preferences $V_1 ,V_2 ,V_3$ in $P$
such that $x_i$ lies between the other two candidates in $V_i$, for all $i=1,2,3$ (\emph{medium-restricted});
\item there do not exist two preferences $V_1 ,V_2$ and four candidates $a,b,c,d$ such that
$a>_{V_1}b>_{V_1}c>_{V_1}d$ and $b>_{V_2}d>_{V_2}a>_{V_2}c$.
\end{itemize}
\end{prop}

Recalling the first remark at the end of the previous section,
since the second condition above is given in terms of a forbidden (4,2)-configuration, which is $(abcd,bdac)$,
we can represent it by means of a set of permutations closed under inversion, which is the set $\{ 2413,3142\}$.
In fact, the permutation mapping the preference $abcd$ into the preference $bdac$ is 2413,
whereas the permutation mapping the preference $bdac$ into the preference $abcd$ is 3124.

\bigskip

We say that $(C,P)$ is an \emph{enriched group-separable election} whenever it is a group-separable election
which, in addition, avoids the (4,2)-configurations $(abcd,badc)$ and $(acbd,bdac)$.
Therefore an enriched group-separable election is a medium-restricted election which further avoids a certain set of (4,2)-configurations.
Such a set can be interpreted as the set of all pairs of linear extensions of the following two posets:

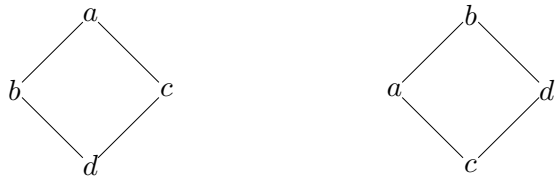
\begin{figure}[h!]
\begin{center}
\begin{tikzpicture}

\begin{scope}

\draw (1,2) node {$a$};
\draw (0,1) node {$b$};
\draw (2,1) node {$c$};
\draw (1,0) node {$d$};

\draw[very thin] (0.9,1.9) -- (0.1,1.1);
\draw[very thin] (1.1,1.9) -- (1.9,1.1);
\draw[very thin] (0.1,0.9) -- (0.9,0.1);
\draw[very thin] (1.9,0.9) -- (1.1,0.1);

\end{scope}

\begin{scope}

\draw (6,2) node {$b$};
\draw (5,1) node {$a$};
\draw (7,1) node {$d$};
\draw (6,0) node {$c$};

\draw[very thin] (5.9,1.9) -- (5.1,1.1);
\draw[very thin] (6.1,1.9) -- (6.9,1.1);
\draw[very thin] (5.1,0.9) -- (5.9,0.1);
\draw[very thin] (6.9,0.9) -- (6.1,0.1);

\end{scope}
\end{tikzpicture}
\end{center}
\end{figure}

More specifically, the four forbidden configurations are $(abcd,badc)$,$(abcd,bdac)$,$(acbd,badc)$ and $(acbd,bdac)$.
In terms of permutations, the above set of forbidden configurations can be represented by the set of permutations
$\Gamma =\{ 2413,3142,2143,3412\}$, which is of course closed under inversion.

\bigskip

The combinatorics of the set $Av(\Gamma )$ of all permutations which avoid the patterns in $\Gamma$ is known and quite nice and interesting.
First of all, the permutations of $Av(\Gamma )$ are precisely those permutations which can be drawn on an X \cite{W}.
The generating function of such permutations (with respect to the length) is $\frac{1-3x}{1-4x+2x^2}$,
and a recurrence relation for the resulting integer sequence is $f_n =4f_{n-1}-2f_{n-2}$, with initial conditions $f_0 =1$ and $f_1 =1$.
From this linear recurrence it is not hard to find a closed formula, which is $f_n =\frac{1}{2}(2+\sqrt{2})^n +\frac{1}{2}(2-\sqrt{2})^n$.
This is a rather interesting sequence from a combinatorial point of view, as it is witnessed by a few papers studying it
(see entry A006012 in \cite{S}).

\bigskip

The elections avoiding the four above configurations (without the medium-restricted condition) can be characterized quite easily.
For any subset $C'\subset C$ of candidates, given an elector $\gamma$, we denote with $E_\gamma ^{C'}$ the 2-element set consisting of
the most preferred and the least preferred candidates of $\gamma$ when restricted to $C'$. Similarly, we denote with
$M_\gamma ^{C'}$ the complement of $E_\gamma ^{C'}$ in $C'$, that is $M_\gamma ^{C'}=C'\setminus E_\gamma ^{C'}$.
The proof of the following proposition is left to the reader.

\begin{prop}
An $(m,n)$-election $(C,P)$ avoids the four above configurations if and only if for every $C'\subseteq C$, with $|C'|=4$,
and for any two electors $\gamma ,\delta$, we have that $E_\gamma ^{C'}\neq M_\delta ^{C'}$.
\end{prop}

Our next aim is to provide a recursive characterization of enriched group-separable elections.

\begin{teor}\label{characterization}
Let $(C,P)$ be an $(m,n)$-election, with $P=(V_1 ,\ldots ,V_n )$.
W.l.o.g., suppose that the preference $V_1$ is represented by the identity permutation of length $m$, that is $V_1 =12\cdots m$.
Then $(C,P)$ is enriched group-separable if and only if there exists $k<m$ such that the multiset of the preferences in $P$ can be partitioned into two blocks as follows:
\begin{itemize}
\item the first block contains all preferences that are isomorphic either to the identity $12\cdots m$ or to its reverse $m\cdots 21$;
\item for the second block, we have two distinct possibilities, which are mutually exclusive:
\begin{itemize}
\item the preferences in the second block are of the form $12\cdots k\pi$ or of the form $\sigma k\cdots 21$,
where $\pi$ and $\sigma$ are permutations of the set $C'=\{ k+1,\ldots ,m\}$ and the restriction of $(C,P)$ to $C'$ is
enriched group-separable, or else
\item the preferences in the second block are of the form $\tau (k+1)\cdots m$ or of the form $m\cdots (k+1)\rho$,
where $\tau$ and $\rho$ are permutations of the set $C'=\{ 1,\ldots k\}$ and the restriction of $(C,P)$ to $C'$ is
enriched group-separable.
\end{itemize}
\end{itemize}
\end{teor}

\emph{Proof.}\quad Suppose first that $(C,P)$ is an enriched group-separable $(m,n)$-election.
From the group-separability we have that there exists a partition of the set of candidates $C$ into two blocks $A$ and $B$ such that
each elector $i$ either ranks $A$ above $B$ or vice versa.
To simplify notations, we will denote these two situations by $A>_i B$ and $A<_i B$, respectively.
Since $V_1 =12\cdots m$, this means that there exists a positive integer $k<m$ such that $A=\{ 1,2,\ldots ,k\}$ and $B=\{ k+1,\ldots ,m\}$.
We start by observing that
the presence in $(C,P)$ of the preferences $12\cdots m$ and $m\cdots 21$ does not prevent the election from being enriched group-separable.
Therefore we can temporarily forget all the preferences of this form and consider all the remaining preferences.

Let $V_i$ be a ranking for which $A>_i B$. Since the permutation generated by the pair $(V_1 ,V_i )$ must avoid $2143$,
$V_i$ cannot contain both an inversion inside its prefix of length $k$ and an inversion inside its suffix of length $m-k$.
Therefore, necessarily in $V_i$ either the elements of $A$ are ranked as in $V_1$ or the elements of $B$ are ranked as in $V_1$.
In other words, either $V_i =1\cdots k\pi$ or $V_i =\tau (k+1)\cdots m$, for suitable permutations $\pi$ and $\tau$.
In a completely analogous way it can be shown that, in case $A<_i B$ in $V_i$, due to the presence of the forbidden pattern $3412$,
$V_i$ cannot contain both a noninversion inside its prefix of length $m-k$ and a noninversion inside its suffix of length $m-k$.
Thus either $V_i =m\cdots (k+1)\rho$ or $V_i =\sigma k\cdots 1$, for suitable permutations $\rho$ and $\sigma$.

Now suppose that $V_i ,V_j$ are two preferences such that $V_i =1\cdots k\pi$ and $V_j =m\cdots (k+1)\rho$, for suitable $\pi$ and $\rho$ as above.
In particular, $\pi$ must contain at least one inversion and $\rho$ must contain at least one noninversion,
otherwise $V_i$ would be the identity and $V_j$ would be the reverse identity.
Therefore the permutation associated with the pair $(V_i ,V_j)$ would contain an occurrence of the pattern $3412$, which is forbidden.
In a similar way, we cannot have two preferences $V_i ,V_j$ such that $V_i =\tau (k+1)\cdots m$ and $V_j =\sigma k\cdots 21$,
otherwise an occurrence of $3412$ would show up as well
(again because $\tau$ must contain at least one inversion and $\sigma$ must contain at least one noninversion).

To conclude the first part of the proof we just observe that, if we restrict our election to any subset of the candidates,
then necessarily the permutations associated with any two preferences cannot contain any pattern of $\Gamma$,
otherwise the whole election would contain such a pattern (this is actually a general property of \emph{configuration definable} elections,
which can be characterized as those elections which are \emph{hereditary}, see \cite{LL} for details).
This is equivalent to saying that the restriction of $(C,P)$ to any subset of candidates is enriched group-separable.

\bigskip

To prove the converse, suppose that the election $(C,P)$ is of the form specified in the statement of the theorem.
In particular, we have two possible types of elections that are allowed.
\begin{enumerate}
\item The multiset of all preferences can be partitioned into the following four subsets:
\begin{itemize}
\item preferences isomorphic to the identity $12\cdots m$ (\emph{type $I$});
\item preferences isomorphic to the reverse identity $m\cdots 21$ (\emph{type $J$});
\item preferences of the form $12\cdots k\pi$ (\emph{type $\alpha$});
\item preferences of the form $\sigma k\cdots 21$ (\emph{type $\delta$});
\end{itemize}
moreover, the restriction of $(C,P)$ to the set of candidates $\{ k+1,\ldots m\}$ is enriched group-separable.
\item The multiset of all preferences can be partitioned into the following four subsets:
\begin{itemize}
\item preferences isomorphic to the identity $12\cdots m$ (\emph{type $I$});
\item preferences isomorphic to the reverse identity $m\cdots 21$ (\emph{type $J$});
\item preferences of the form $\tau(k+1)\cdots m$ (\emph{type $\beta$});
\item preferences of the form $m\cdots (k+1)\rho$ (\emph{type $\gamma$});
\end{itemize}
moreover, the restriction of $(C,P)$ to the set of candidates $\{ 1,\ldots k\}$ is enriched group-separable.
\end{enumerate}

We now show in some detail that the elections of the first type are enriched group-separable.
The proof for elections of the second type can be done in a completely analogous way.

Since it is clear that, setting $A=\{ 1,2,\ldots ,k\}$ and $B=\{ k+1,\ldots ,m\}$, for any preference $V_i$ in $P$ we have either $A<_i B$ or $B<_i A$,
the election $(C,P)$ is group-separable.
To conclude the proof, it will be enough to show that $(C,P)$ avoids the two configurations represented by the permutations 2413 and 3412.

Let $V_i ,V_j$ be two preferences in $P$.
Given four candidates in the set $\{ k+1,\ldots ,m\}$, they cannot constitute a forbidden configuration,
since by hypothesis the restriction of the election to $\{ k+1,\ldots ,m\}$ is enriched group-separable.
So suppose that $x,y,z,t$ are any four candidates, with $x\in \{ 1,2, \ldots ,k\}$.
A case by case analysis shows that $x,y,z,t$ cannot constitute a forbidden configuration, whatever the types of $V_i$ and $V_j$ are.
More specifically,
they always form a configuration which must have one of the two elements 1 and 4 either at the beginning or at the end
(and so cannot be neither 2413 nor 3412).
For instance, if $V_i$ is of type $J$ and $V_j$ is of type $\alpha$,
then the smallest among $x,y,z,t$ is the rightmost of them in $V_i$ and the leftmost of them in $V_j$
(this is due to the fact that there is at least one of such elements which is $\leq k$).
Therefore in the pair $(V_i ,V_j)$ the four elements $\{ x,y,z,t\}$ constitute a configuration starting with 4.
All the remaining cases can be analyzed in a completely analogous way.\cvd

The previous theorem provides an interesting characterization of enriched group-separable elections.
First of all, an enriched group-separable election is a special group-separable election,
and so there are two blocks of candidates $A$ and $B$ such that each elector ranks either $A$ above $B$ or vice versa.
This means that the set of electors essentially splits into two groups,
whose respective preferences go either to the candidates of $A$ or to the candidates of $B$.
Moreover, there is an additional requirement for enriched elections, which is the following (possibly exchanging the roles of $A$ and $B$):
all the electors which prefer $A$ to $B$ rank the candidates of $A$ in exactly the same way;
moreover, all the remaining electors (i.e., those who prefer $B$ to $A$) rank the candidates of $A$ precisely in the reverse order.

The special form of enriched group-separable elections allows to find a recurrence relation for their enumeration, from which a closed formula can be deduced.

\begin{prop}
Denote with $f(m,n)$ the number of enriched group-separable elections with $m$ candidates and $n$ electors
and with $f_r (m,n)$ the number of such elections in which one of the preferences is fixed.
Then the following recurrence relation holds:
\begin{equation}\label{recurrence}
f_r (m,n)=2^n f_r (m-1,n)-2^{n-1}f_r (m-2,n),
\end{equation}
with initial conditions $f_r (0,n)=f_r (1,n)=1$.
Moreover $f(m,n)=m!\cdot f_r (m,n)$.
\end{prop}

\emph{Proof.}\quad The fact that $f(m,n)=m!\cdot f_r (m,n)$ is due to the obvious fact that
the total number of enriched group-separable elections is given by the number of enriched group-separable elections in which one preference is fixed
times the number of possible choices for that preference.

Now suppose, w.l.o.g., that the first preference is fixed and equal to $12\cdots m$.
Call \emph{reduced} such an enriched group-separable election and denote it with $(C,P)$.
There are two distinct possible types for $(C,P)$, for each of which the multiset of preferences can be partitioned into two classes.
Specifically, the two types are the following:
\begin{itemize}
\item either a preference is of the form $1\pi$ or of the form $\rho 1$, with $\pi ,\rho$ of length $m-1$,
and $(C,P)$ restricted to the set of candidates $\{ 1,2,\ldots m-1\}$ is enriched group-separable (and, of course, reduced);
\item either a preference is of the form $\sigma m$ or of the form $m\tau$, with $\sigma ,\tau$ of length $m-1$,
and $(C,P)$ restricted to the set of candidates $\{ 2,3,\ldots m\}$ is enriched group-separable (and, of course, reduced).
\end{itemize}

The above statement is a consequence of Theorem \ref{characterization}.
In particular, it is not difficult to realize that, depending on which of the two above types $(C,P)$ is,
removing either $m$ or $1$ results in an enriched group-separable election.

For both of the above types, in order to determine $(C,P)$,
one has just to decide which subset of electors is associated with preferences of one of the two forms.
Therefore, after having chosen one of the $2^n$ possible subsets of electors,
we simply have to choose a (reduced) enriched group-separable election on a subset of candidates of cardinality $m-1$.
This gives a total of $2^n f_r (m-1,n)$ elections. However, the two above types of elections are not disjoint.
Indeed, elections all of whose preferences have the form $1\psi m$ are of both types, so they are counted twice.
The total number of such elections is clearly $2^{n-1}f_r (m-2,n)$,
since it is obvious that the restriction to the set of candidates $\{ 2,3,\ldots ,m-2,m-1\}$ is (reduced) enriched group-separable.
Subtracting such a quantity we obtain precisely formula (\ref{recurrence}), as desired.\cvd

Equation (\ref{recurrence}) is a linear recurrence relation in the indeterminate $m$ with constant coefficients.
It can be solved using standard methods, giving the following closed expression for its general term:
$$
f_r (m,n)=\frac{\varphi+(1-2^{n-1})}{2\varphi}\left( 2^{n-1}+\varphi \right)^m +\frac{\varphi-(1-2^{n-1})}{2\varphi}\left( 2^{n-1}-\varphi \right)^m ,
$$
where $\varphi=\sqrt{2^{n-1}(2^{n-1}-1)}$.
Since $f(m,n)=m!\cdot f_r (m,n)$,
we thus have a closed form for the number $f(m,n)$ of enriched group-separable elections with $m$ candidates and $n$ electors.

For small values of $m$ and/or $n$, we are able to write some interesting explicit formulas. The following corollary is immediate, so the proof is omitted.

\begin{cor}
\begin{enumerate}
\item $f(3,n)=6\cdot 2^{n-1}(2^n -1)$.
\item $f(4,n)=24\cdot 4^{n-1}(2^{n+1}-3)$.
\item $f(5,n)=120\cdot 4^{n-1}(2^{2n+1}-2^{n+2}+1)$.
\item $f(m,2)=\frac{m!}{4}\cdot\left( (2+\sqrt{2})(2-\sqrt{2})^m +(2-\sqrt{2})(2+\sqrt{2})^m \right)$.
\end{enumerate}
\end{cor}

At least two of the above formulas are worth a comment.

Rather surprisingly,
the expression we have found for $f(4,n)$ coincides with the number of single-peaked elections with the same number of candidates and electors.
We do not have any explanation of this fact, it would be interesting to understand if there is a connection between the two types of elections.
Notice however that the enriched group-separable domain is in general different from the single-peaked domain.
This can be shown, for instance, by counting elections of the two types when $m=5$ and $n=2$.
In fact, the results from \cite{LL} imply that, in such a case, there are 8400 single-peaked elections,
whereas the above corollary tells us that there are 8160 enriched group-separable elections.

Concerning the formula for $f(m,2)$, observe that when $n=2$ the medium-restricted configuration cannot appear,
so the enumeration of the associated elections coincides with the enumeration of permutations avoiding the patterns in $\Gamma$.
Indeed, recurrence (\ref{recurrence}) becomes $f(m,2)=4f(m-1,2)-2f(m-2,2)$, the same recurrence we have mentioned at the beginning of this section.

\section{Configuration containment and strong order on pairs of permutations}

Given permutations $\pi ,\rho \in S_n$ and $\tau ,\sigma \in S_h$, with $h\leq n$,
we say that $[\tau ,\sigma ]\leq [\pi ,\rho ]$ whenever
there exist $h$ distinct values $\alpha_1 ,\ldots ,\alpha_h \leq n$
which constitute an occurrence of both $\tau$ in $\pi$ and $\sigma$ in $\rho$.

For instance, $[213,132]\leq [614235,126534]$,
since the elements 2,4 and 5 form an occurrence of 213 in 614235
and the same elements form an occurrence of 132 in 126534.
On the other hand, $[213,132]\nleq [614235,152436]$.
In this case, we say that $[614235,152436]$ \emph{avoids} $[213,132]$.

The above relation is clearly a partial order on pairs of permutations (i.e., on the set $S^2$)
which we call the \emph{strong order}.
This is motivated by the fact that this partial order is of course stronger than the usual componentwise order
of the direct product of the pattern poset with itself.

The strong order on pairs of permutations is not a new notion.
It was first defined in \cite{AADHHM} under the name of involvement order (on pairs of permutations),
and some of its properties has also been investigated.
However, there is some evidence that this partial order is not an easy one,
and in particular is much harder than the analogous pattern order on permutations.
For instance, the problem of enumerating classes of pattern-avoiding pairs of permutations is largely open
also for very simple patterns.
Just to mention the simplest non trivial case, it is easy to realize that there is only one pattern of length 2
up to symmetries, which is $(12,21)$.
However, given a pair of permutations $(\pi ,\rho )$ of length $n$,
it can be shown that $(\pi ,\rho )$ avoids $(12,21)$ if and only if $\rho$ is below $\pi$ in the weak Bruhat order:
the problem of counting such pairs is a well known difficult problem, which is still unsolved.

It turns out that the notion of strong order on pairs of permutations shows up in a natural way to express
containment of configurations in elections.

\bigskip

The next theorem generalizes the result obtained by Lackner and Lackner (recorded here as Proposition \ref{LandL}) to the case of elections with 3 voters.

\begin{teor}\label{3-election}
Let $(C,P)$ be a $(m,3)$-election, with $C=\{ x_1 ,\ldots ,x_m \}$ and $P=(V_1 ,V_2 ,V_3 )$.
Let $(S,T)$ be a $(h,3)$-configuration, with $S=\{ s_1 ,\ldots ,s_h \}$ and $T=(T_1 ,T_2 ,T_3 )$.
Moreover, suppose that the preferences are given as follows:
\begin{align*}
V_1 &: x_1 x_2 \cdots x_m & \qquad V_2 &: x_{\pi (1)} x_{\pi (2)} \cdots x_{\pi (m)} & \qquad
V_3 &: x_{\rho (1)} x_{\rho (2)} \cdots x_{\rho (m)} \\
T_1 &: s_1 s_2 \cdots s_h & \qquad T_2 &: s_{\tau (1)} s_{\tau (2)} \cdots s_{\tau (h)} & \qquad
T_3 &: s_{\sigma (1)} s_{\sigma (2)} \cdots s_{\sigma (h)},
\end{align*}
for certain permutations $\pi ,\rho ,\tau ,\sigma$. Then $(C,P)$ contains $(S,T)$ if and only if
$[\pi ,\rho ]$ contains one of the following patterns:
$$
[\tau, \sigma ], [\sigma ,\tau ], [\tau^{-1},\tau^{-1}\circ \sigma ], [\tau^{-1}\circ \sigma ,\tau^{-1}],
[\sigma^{-1},\sigma^{-1}\circ \tau ], [\sigma^{-1}\circ \tau ,\sigma^{-1}].
$$
\end{teor}

\emph{Proof.} For any element $x\in C$, we will say that the \emph{index} of $x$ is $i$ when $x=x_i$.
A similar terminology will be used for the elements of $S=\{ s_1,\ldots ,s_h \}$.

\begin{itemize}

\item[$\Leftarrow$)] Suppose that $[\pi ,\rho ]$ contains $[\tau, \sigma ]$.
This means that there exist elements $\alpha_1 ,\alpha_2 ,\ldots ,\alpha_h \in C$
whose indices form an occurrence of $\tau$ into $\pi$ and of $\sigma$ into $\rho$.
It is not restrictive to suppose that
the $\alpha_i$'s are listed in increasing order of their indices (as elements of $C$),
so that, if $i<j$, then $\alpha_i <_{V_1}\alpha_j$, for all $i,j\leq h$.
Therefore, choosing $f:\{ 1,2,3\}\rightarrow \{ 1,2,3\}$ as the identity function and $g:S\rightarrow C$ such that $g(s_i )=\alpha_i$,
it is easy to see that $(S,T)$ is contained in $(C,P)$.
Indeed, if $s_i <_{T_1}s_j$, then $i<j$, and so $g(s_i )=\alpha_i <_{V_1}\alpha_j =g(s_j )$;
moreover, if $s_{\tau (i)}<_{T_2}s_{\tau (j)}$, then $g(s_{\tau (i)})<_{V_2}g(s_{\tau (j)})$
because by hypothesis the indices of the $\alpha_i$'s form an occurrence of $\tau$ in $\pi$;
for a similar reason we also have that,
if $s_{\sigma (i)}<_{T_3}s_{\sigma (j)}$, then $g(s_{\sigma (i)})<_{V_3}g(s_{\sigma (j)})$.

If instead $[\pi ,\rho ]$ contains $[\sigma, \tau]$, we can use the same argument as above,
just interchanging the roles of $\sigma$ and $\tau$.
More formally, we just have to choose as $f$ the function mapping 1 into 2, 2 into 1 and 3 into 3.

Now suppose that $[\pi ,\rho ]$ contains $[\tau^{-1},\tau^{-1}\circ \sigma ]$.
Again, this means that there exist elements $\alpha_1 ,\alpha_2 ,\ldots ,\alpha_h \in C$
whose indices form an occurrence of $\tau^{-1}$ into $\pi$ and of $\tau^{-1}\circ \sigma$ into $\rho$.
As before, we can assume that, if $i<j$, then $\alpha_i <_{V_1}\alpha_j$, for all $i,j\leq h$.
Let $f:\{ 1,2,3\}\rightarrow \{ 1,2,3\}$ be the function mapping 1 into 2, 2 into 1 and 3 into 3;
moreover, define $g:S\rightarrow C$ by setting $g(s_i )=\alpha_{\tau^{-1}(i)}$.
First of all we observe that, if $s_i <_{T_1}s_j$, then $i<j$,
and so $g(s_i )=\alpha_{\tau^{-1}(i)}<_{V_2}\alpha_{\tau^{-1}(j)}=g(s_j )$,
since by hypothesis the indices of the $\alpha_i$'s form an occurrence of $\tau^{-1}$ in $\pi$.
Moreover, if $s_{\tau(i)}<_{T_2}s_{\tau(j)}$, then $g(s_{\tau(i)})=\alpha_i <_{V_1}\alpha_j =g(s_{\tau(j)})$,
as a consequence of our assumption on the $\alpha_i$'s.
Finally, if $s_{\sigma(i)}<_{T_3}s_{\sigma(j)}$,
then $g(s_{\sigma(i)})=\alpha_{\tau^{-1}(\sigma(i))}<_{V_3}\alpha_{\tau^{-1}(\sigma(j))}=g(s_{\sigma(j)})$,
since the indices of the $\alpha_i$'s form an occurrence of $\tau^{-1}\circ \sigma$ in $\rho$.
We have thus shown that, in this situation, $(C,P)$ contains $(S,T)$.

The analysis of all the remaining cases can be obtained by the previous ones
by suitably interchanging the role of the permutations involved and choosing the appropriate function $f$,
so it is left to the reader.

\item[$\Rightarrow$)] Suppose that $(C,P)$ contains $(S,T)$.
This means that there are suitable functions $f$ and $g$
that satisfy the definition of configuration containment given in the Introduction.
Depending on the particular choice of $f$ and $g$, there are several cases to analyze.
We will provide details for only two of them,
since, as before, all the remaining ones can be easily derived from those.

Suppose first that $f$ is the identity map on $\{ 1,2,3\}$ and set $g(s_i )=\alpha_i$ , for all $i\leq h$.
In this case, the hypothesis means exactly that the indices of the $\alpha_i$'s form an occurrence of $\tau$ in $\pi$ and of $\sigma$ in $\rho$,
i.e $[\tau ,\sigma ]\leq [\pi ,\rho ]$.

Now suppose that $f$ maps 1 into 2, 2 into 1 and 3 into 3, and set $g(s_i )=\alpha_i$ , for all $i\leq h$.
We wish to show that $[\tau^{-1},\tau^{-1}\circ \sigma ]\leq [\pi ,\rho ]$.
First of all, since $f$ maps 2 into 1, we have that,
in the sequence $\alpha_{\tau (1)},\alpha_{\tau (2)},\ldots ,\alpha_{\tau (h)}$,
elements are listed in increasing order of their indices (as elements of $C$).
This means that,
denoting again by $\alpha$ the permutation isomorphic to the word defined by the indices of the word
$\alpha_1 \cdots \alpha_h$ (with a slight abuse of notation),
the composition $\alpha \circ \tau$ is (isomorphic to) the identity permutation on the set of the indices of the $x_i$'s.
Equivalently, the word $\alpha_1 \cdots \alpha_h$ determines a permutation isomorphic to $\tau^{-1}$ on the set of the indices of the $x_i$'s appearing in the set of the $\alpha_i$'s.
Thus, since $f$ maps 1 into 2, we get that, looking at the indices, $\alpha_1 \cdots \alpha_h$ shows an occurrence of $\tau^{-1}$ in $\pi$,
which was the first thing to prove.
Moreover, since $f$ maps 3 into 3, we have that $x_{\rho (1)} x_{\rho (2)} \cdots x_{\rho (m)}$ contains
$\alpha_{\sigma (1)} \alpha_{\sigma (2)} \cdots \alpha_{\sigma (h)}$ as a subword
(whose letters are not necessarily contiguous, of course).
Since, as we have seen before,
the word $\alpha_1 \cdots \alpha_h$ induces on the indices the permutation $\tau^{-1}$,
looking at the indices of the $\alpha_i$'s inside $x_{\rho (1)} x_{\rho (2)} \cdots x_{\rho (m)}$,
we observe that they form an occurrence of $\tau^{-1}\circ \sigma$ inside $\rho$, which is what remained to be proved.\cvd

\end{itemize}

We now state two corollaries of the previous theorem. The first one gives a special case, which appears to be more manageable.
The second one suggests an approach for finding bounds to the number of elections defined in terms of forbidden configurations,
such that at least one of the configurations is an election with 3 voters.

\begin{cor}
In the hypothesis of the above proposition, suppose that $\tau$ is the identity permutation.
Then $(C,P)$ contains $(S,T)$ if and only if $[\pi ,\rho ]$ contains one of the following patterns:
$$
[id,\sigma ],[\sigma ,id],[\sigma^{-1},\sigma^{-1}],
$$
where $id$ denotes an identity permutation (of the appropriate length).
\end{cor}

\begin{cor}\label{3-conf}
Let $V_1$ be a total order on $C$ (with $|C|=m$). Let $(S,T)$ be a $(h,3)$-configuration as in Theorem \ref{3-election}.
The number of pairs $(V_2 ,V_3 )$ of total orders on $C$ such that $(C,P)$ avoids $(S,T)$ (with $P=(V_1 ,V_2 ,V_3 )$) is equal to
$|S_m (\Pi )|$, where $\Pi$ is the set of pairs of permutations listed in Theorem \ref{3-election}.
\end{cor}

The last corollary generalizes Corollary \ref{inverse}, which has to do with $(h,2)$-configurations, to the case of $(h,3)$-configurations.
At least in principle, our result could allow to find some upper bound for elections avoiding at least one $(h,3)$-configuration.
Suppose indeed that $a(m,n,\Gamma)$ denotes the number of $(m,n)$-elections $(C,P)$, with $P=(V_1 ,\ldots V_n)$ as usual,
avoiding a set of configurations $\Gamma$ and let $(S,T)\in \Gamma$ be a $(h,3)$-configuration.
In order to determine an upper bound for $a(m,n,\Gamma)$, we can find an upper bound for the number of $(m,n)$-elections avoiding $(S,T)$.
First of all, we can choose the first ranking $V_1$ of $(C,P)$ at random, hence we have $m!$  possibilities.
Then we can relax the condition of avoiding $\Gamma$ by simply requiring for $(C,P)$ to avoid the $(h,3)$-configuration $(S,T)$.
Even more specifically, we ask that for none of the triples $(V_1 ,V_i ,V_j)$, with $1<i<j$, the election $(C,(V_1 ,V_i ,V_j ))$ contains $(S,T)$.
Using Corollary \ref{3-conf} and its notations, we know that there are $|S_m (\Pi )|$ ways to choose each pair $(V_i ,V_j )$,
where $\Pi$ is the set of pairs of permutations listed in Theorem \ref{3-election}.
Now, since the possible number of pairs $(i,j)$, with $1<i<j$, is ${n-1\choose 2}$, we find the following upper bound:
\begin{equation}\label{upper_bound_3config}
a(m,n,\Gamma)\leq m!|S_m (\Pi )|^{{n-1\choose 2}}.
\end{equation}

To illustrate how to apply the above formula, we sketch how it reads in the case of the single-crossing domain.
Recall that preferences are \emph{single-crossing} if there exists a linear ordering of the voters such that,
for any pair of candidates along this ordering, there is a single spot where the voters switch from preferring one candidate above the other one.
In \cite{BCW} single-crossing domains are characterized in terms of two forbidden configurations, having three voters and four voters respectively.
Since there is no forbidden configuration with 2 voters, Corollary \ref{inverse} of \cite{LL} is not useful here.
On the other hand, since one of the forbidden configurations has three voters, we can describe an upper bound using formula (\ref{upper_bound_3config}).
Without delving into details
(in the forbidden configurations of the characterization of \cite{BCW} some of the candidates are not necessarily distinct),
we can show that single-crossing elections must avoid a certain configuration with three voters which, in the notations of Corollary \ref{3-conf},
is associated with the set of permutations
$$
\Pi=\{ [4231,4132],[4132,4231],[4231,1432],[1432,4231],[2431,1432],[1432,2431]\} .
$$

Thus we have that an upper bound for the number of single-crossing $(m,n)$-elections is given by formula (\ref{upper_bound_3config}),
where $\Pi$ is the set of permutations displayed above.

Unfortunately, as we recalled at the beginning of the present section,
the enumerative combinatorics of pattern avoidance for the strong order on pairs of permutations is still largely unknown.
We hope that our work, showing the connections with elections defined by forbidden configurations, can stimulate further research on such a topic.

\section{Conclusions}

Following the footprints of \cite{LL},
we have exploited the connections between permutation patterns and forbidden configurations in restricted elections
in order to further show how fruitful such a connection can be.
Our first result has been the definition of a new type of domain restriction, called \emph{enriched group-separable},
for which we have been able to provide complete enumerative results
(in an easier way than it has been done for the more general group-separable restriction).
In the same way, it is conceivable that, starting from suitable sets of permutations closed under taking the inverse,
further natural and interesting domain restrictions can be defined, for which the exact enumeration could be possible.
Our second result has been a generalization of the main result of \cite{LL}
concerning the above mentioned connections between permutation patterns and restricted elections.
Our generalization involves the notion of strong containment for pairs of permutations,
and suggests a possible approach for finding upper bounds on the number of restricted elections
when the defining forbidden configurations contains at least one configuration with 3 votes.
This improves the approach of \cite{LL}, which works only when there is a forbidden configuration with 2 voters.
Unfortunately, our approach requires a deeper knowledge of the enumerative combinatorics of pattern avoiding pairs of permutations,
which is presently unavailable. We thus hope that our achievement will encourage further research in this direction.


\end{document}